\def\nubar{\bar {\nu}}
\def\dbar{\bar {d}}
\def\etabar{\bar {\eta}}
\documentstyle[amssymb]{jtbart}
\newbox\smilebox
\newbox\anchorbox
\newbox\noanchorbox
\newbox\tempbox

\setbox\smilebox=\hbox{$\smile$}

\def\anchor{\hbox{\vtop{
           \hbox to \wd\smilebox{\hfil\vrule width.4pt height7pt depth1pt\hfil}
           \vskip  -11.5truept
           \hbox to \wd\smilebox{\hfil$\smile$\hfil}}}}
\setbox\anchorbox=\anchor
\def\noanchor{\hbox{\vtop{
           \hbox to \wd\anchorbox{\hfil\anchor\hfil}
           \vskip -14truept
           \hbox to \wd\anchorbox{\hfil/\hfil}}}}
\setbox\noanchorbox=\noanchor

\def\fg#1#2#3{\setbox\tempbox=\hbox{$\scriptstyle{#2}$}
\ifnum\wd\anchorbox>\wd\tempbox\dimen255=\wd\anchorbox
\else\dimen255=\wd\tempbox\fi
{#1\,\vtop{\hbox to \dimen255{\hfil\anchor\hfil}
           \vskip -6truept
           \hbox to \dimen255{\hfil$\scriptstyle{#2}$\hfil}}
           \,#3}}

\def\bet{]}

\def\lg{\rm lg }

\def\lev{\rm lev }

\def\Aun{\under {A}}
\def\conc{\widehat{~~}}

\makeatletter
\def\@begintheorem#1#2{\it \trivlist \item[\hskip
\labelsep{\bf #2\ #1.}]}
\def\@opargbegintheorem#1#2#3{\it \trivlist
      \item[\hskip \labelsep{\bf #1\ #2\ (#3).}]}
\makeatother

\newtheorem{them}{Theorem}[section]

\newtheorem{lemm}[them]{Lemma}

\newtheorem{prop}[them]{Proposition}
\newtheorem{thm}[them]{Theorem}

\newcommand{\jtbnumpar}[1]{\refstepcounter{them}
\trivlist
\item[\hskip \labelsep{\bf \thethm \ #1.}]}

\newcommand{\jtbdef}{\jtbnumpar{Definition}}
\def\jtbnot{\jtbnumpar{Notation}}

\def\PROOF #1.{\par\noindent{\it Proof#1}.\ \ignorespaces}
\def\proof #1.{\par\noindent{\it Proof#1}.\ \ignorespaces}

\def\Mscr{{\cal M}}

\let\?=\joinrel

\let\meet=\bigwedge

\let\leftv=^
\let\rightv=^

\def\Pscr{{\cal P}}

\def\k/{\kern.2em}    

\def\dom{\mathop{\rm dom}}

\def\acl{\mathop{\rm acl}}

\def\sup{\mathop{\rm sup}}

\def\dim{\mathop{\rm dim}}

\def\lg{\mathop{\rm lg}}

\def\stp{\mathop{\it stp}}

\def\acl{{\rm acl}}     

\let\bet=\Im        %

 \mathchardef\bet="0B69

\let\union=\cup             %

\edef\bigcup{\mathop{\textstyle\mathchar\the\bigcup}}
\let\bigunion =\bigcup

\let\inter=\cap             %

\edef\bigcap{\mathop{\textstyle\mathchar\the\bigcap}}

\edef\bigwedge{\mathop{\textstyle\mathchar\the\bigwedge}}

\edef\bigvee{\mathop{\textstyle\mathchar\the\bigvee}}

\edef\sum{\mathop{\textstyle\mathchar\the\sum}}

\def\nubar{\overline \nu}

\def\etabar{\overline \eta}
\def\math&{\ \& \ }

\def\force {\mathrel^\joinrel\rightarrow}
\def\force {\mathrel{\scriptstyle\mathrel^\joinrel\rightarrow}}
\def\forceq {\mathrel{\mathop{\force}\limits_{\textstyle\texsim}}}
\def\forceq{\mathrel^\joinrel
 \mathrel{\mathop{\rightarrow}\limits_{\smash{\textstyle\texsim}}}}
\def\forceq{\mathrel{\scriptstyle\mathrel^\joinrel
 \mathrel{\smash{\mathop{\rightarrow}\limits_{\smash{\raise
 2pt\hbox{$\scriptstyle\texsim$}}}}}}}
\let\iff\leftrightarrow         %

\let\exclaim=!                 %
\let\iso\approx             %
\let\texsim=\sim         %
\let\conj\sim             %
\def\conjp #1 {\conj_{#1}}     %
\let\sim\simeq             %
\def\con{{^\frown}}
\let\neg=\lnot             %

\def\0bar{\bar 0}         %
\def\1bar{\bar 1}         %
\def\abar{\overline a}

\let\sat=\models                %

\def\bbar{\overline b}
\def\cbar{\overline c}
\def\dbar{\overline d}
\def\ebar{\overline e}

\def\hbar{\overline h}

\def\ubar{\overline u}
\def\vbar{\overline v}

\def\xbar{\overline x}
\def\ybar{\overline y}

\def\Bun{\underline B}
\def\Aun{\underline A}

\def\tp{\rm tp}

\def\Proof{Proof}

\def\bet{]}

\def\bet{]}

\def\bigunion{\union}

\def\text#1{\ifmmode\leavevmode\hbox{#1}\else
   \typeout{Warning: \string\text \space used outside math mode!}
   \begingroup\hbox{#1}\endgroup\fi}

\author{
J. Baldwin
\thanks{Partially supported by N.S.F. grant 90000139.}
\\Department of Mathematics\\
University of Illinois, Chicago
\and M. C. Laskowski\thanks{Visiting U.I.C. from the University of
Maryland thanks to the NSF Postdoctoral program.}
\\Department of Mathematics\\
University of Illinois, Chicago
\and
S. Shelah\\
Department of Mathematics\\
Hebrew University of Jerusalem
\thanks{The authors thank the U.S. Israel
Binational Science Foundation for its support of this project.
This is item 464 in Shelah's bibliography.}}
\title{Forcing Isomorphism}
\begin{document}
\maketitle
If two models of a first order theory are  isomorphic then they remain
isomorphic in any forcing extension of the universe of sets.  In
general, however, such a forcing extension may create new isomorphisms.
For example, any forcing that collapses cardinals may easily make
formerly
non-isomorphic models isomorphic.  Certain model theoretic constraints
on the theory and other constraints on the forcing can  prevent this
pathology.

A countable
first order theory is said to be {\em classifiable} if it is
superstable and does not have either the dimensional order property
(DOP) or
the omitting types order property (OTOP).  Shelah has shown
\cite{Shelahbook2nd} that if a theory $T$ is classifiable then  each
model of cardinality $\lambda$ is described by a sentence of $L_{\infty,
\lambda}$.  In fact this sentence can be chosen in the
$L^*_{\lambda}$.  ($L^*_{\lambda}$ is the result of enriching the language
$L_{\infty,\beth^+}$ by adding for each $\mu < \lambda$ a
quantifier saying the dimension of a
dependence structure is greater than $\mu$.)  Further work
(\cite{Hartthesis}, \cite{BuechlerShelah})
shows that $\beth^+$ can be replaced by $\aleph_1$.
The truth of such sentences will be preserved by any forcing that
does not collapse cardinals $\le\lambda$ and that
adds no new countable
subsets of $\lambda$,
e.g., a $\lambda$-complete forcing.
That is, if two models of a classifiable theory of power $\lambda$
are non-isomorphic, they remain non-isomorphic after a $\lambda$-complete
forcing.

In this paper we show that
the hypothesis of the forcing
adding no new countable subsets of $\lambda$ cannot be eliminated.
In particular, we show that non-isomorphism of models of a classifiable
theory need
not be preserved by ccc forcings.
The following definition isolates the key issue of this paper.

\jtbdef  Two structures $M$ and $N$ are {\em potentially isomorphic} if
there
is a ccc-notion of forcing $\Pscr$ such that if $G$ is $\Pscr$-generic
then $V[G] \sat M \iso N$.

\smallskip
In the first section we will show that any theory that is not
classifiable has models that are not isomorphic but are potentially
isomorphic.  In the second, we show that this phenomenon
can also
occur for classifiable theories.
The reader may find it useful to examine first the example
discussed in Theorem~\ref{examp}.
\section{Non-classifiable Theories}
\label{unclass}
We begin by describing a class (which we call amenable) of subtrees of
$Q^{\leq\omega}$ that are pairwise potentially isomorphic.   Then we use
this fact to
show that every nonclassifiable theory has a pair of models that are not
isomorphic but are potentially isomorphic.

\jtbnot \mbox{}
\label{maindef}
\begin{enumerate}
\item
We adopt the following notation for relations on subsets of
$Q^{\leq\omega}$.  $\sqsubset$ denotes the subsequence relation;
$<$ denotes lexicographic ordering; for $\alpha \leq \omega$,
$\lev_{\alpha}
$ is a unary predicate that holds of sequences of
length (level) $\alpha$; $\wedge$ is the operation on two
sequences that produces their largest common initial segment.
We denote the ordering of the rationals by $<_Q$.
\item
For $\eta \in Q^{\omega}$, let
$D_{\eta} = \{\sigma \in Q^{\omega}:\sigma(2n) =
\eta(n)\}$ and
$S_{\eta} = \{\sigma \in D_{\eta}:\sigma(2n+1)
\text{ is 0 for all but finitely many $n$}\}$.
Let $C = \bigunion_{\eta\in Q^{\omega}}S_{\eta}$.
\item The language $L^t$ (for tree) contains the symbols $\sqsubset$,
$<$, $\lev_{\alpha}$ and unary predicates $P_{\eta}$ for $\eta \in
Q^{\omega}$.
\item For any $A \subseteq C$, $A^*$ denotes the $L^t$-structure with universe
$A \union Q^{<\omega}$ under the natural interpretation of $\sqsubset, <,
\lev_{\alpha}$ and with $P_{\eta}(A^*) = S_{\eta}\inter A$.
\end{enumerate}

Note that $\langle C, < \rangle$ is isomorphic to a subordering of the
reals. Since $C$
is dense we may assume $Q$ is embedded in $C$ but not necessarily in a
natural way.

\jtbdef  A substructure $A^*$ of $C^*$
is {\em amenable} if for all $\eta \in Q^{\omega}$, all $n\in \omega$
and all $s \in Q^n$, if $P_{\eta}(C^*)$ contains an element extending
$s$ then $P_{\eta}(A^*)$ does also.

\jtbnumpar{Remark}
\label{recastamen}
It is easy to see that a substructure $A^*$ of $C^*$
is {\em amenable} just if
for all even $n$ and all $s \in Q^n$, if $\eta(i) = s(2i)$ for all
$i<{\frac{n}{2}}$ then for
every $r \in Q$
there is a $\nu\in A \inter S_{\eta}$ with $\nu|n+1= s\con r$.

\jtbnumpar{Main Lemma}
\label{ml}
If $A^*$ and $B^*$ are amenable substructures of $C^*$ then they
are
potentially isomorphic.

\Proof.
Let $\Pscr$ denote the set of finite partial $L^t$-isomorphisms between
$A$ and $B$ under the natural partial order of extension.
We can naturally extend any $L^t$ elementary bijection between
$A$ and $B$ to an isomorphism of $A^*$ and $B^*$ as
$L^t$-structures.
\jtbnumpar{Claim 1}  $\Pscr$ satisfies the countable chain condition.

In fact, we will show
$\Pscr = \cup_{n\in \omega} F_n$ where if $p,q \in F_n$ then $p
\union q \in \Pscr$.  Given $p \in \Pscr$, let
$\langle a_1 \ldots a_n \rangle$ be the lexicographic
enumeration
of $\dom p$.  Let $n(p)$ be
the cardinality of $\dom p$ and let $k =k(p)$ be the least integer
satisfying the following conditions.

\begin{enumerate}
\item If $i \neq j$, $a_i|k \neq a_j|k$.
\item If $i \neq j$, $p(a_i)|k \neq p(a_j)|k$.
\item For all $n$ with $2n +1 \geq k$, $a_i(2n+1) = 0$.
\item For all $n$ with $2n +1 \geq k$, $p(a_i)(2n+1) = 0$.
\end{enumerate}

Now define an equivalence relation on $\Pscr$ by $p \sim q$ if $n(p) =
n(q)$, $k(p) = k(q)$ and
letting
$\langle a_1 \ldots a_n \rangle$ enumerate
(in lexicographic order)
$\dom p$,
$\langle a'_1 \ldots a'_n \rangle$ enumerate
(in lexicographic order)
$\dom q$, for each $i$,
$a_i|k(p) = a'_i|k(p)$
and
$p(a_i)|k(p) = q(a'_i)|k(p)$.  Then since $Q$ is countable and the
domains of elements of $\Pscr$ are finite,
$\sim$ has only countably many
equivalence classes; we designate
these classes as the $F_n$.  We must show
that if $p \sim q$ then $p \union q \in \Pscr$.  It suffices to
show that for all $i,j$ if
$C \sat a_i<a'_j$ then
$C \sat p(a_i)<q(a'_j)$.
\begin{itemize}
\item{Case 1:} $i\neq j$.  By the definitions of $k = k(p)$ and
$\sim$, we must have $j > i$,
$a_i|k <
a'_j|k$ and
$p(a_i)|k <
q(a'_j)|k$.   This suffices.
\item{Case 2:} $i = j$.
Choose the least $t$ such that $a_i(t) \neq a'_i(t)$.  Then
$a_i(t) < a'_i(t)$.
Note that $t > k$ and $t$ must be even since for any odd $t> k$,
$a_i(t) = a'_i(t)= 0$.

Suppose $a_i \in S_{\nu}$ and $a'_i \in
S_{\eta}$.

We now claim
$p(a_i)|t = q(a'_i)|t$.  Fix any $\ell < t$.
If $\ell < k$,
$p(a_i)(\ell) = q(a'_i)(\ell)$
by the definition of $\sim$.
If $\ell \geq k$ is odd
then the fourth condition in the
definition of $k(p)$ guarantees that
$p(a_i)(\ell) = q(a'_i)(\ell)= 0$.
Finally, if $\ell \geq k$ and $\ell = 2u$, since $p$ and $q$
preserve the $P_{\eta}$,
$p(a_i)(\ell) = \nu(u)= a_i(\ell)$ and
$q(a'_i)(\ell) = \eta(u)= a'_i(\ell)$
but these are equal by the minimality of $t$.

It remains to show
$p(a_i)(t) <_Q q(a'_i)(t)$.
By condition iv) it follows that $t$ is even.
But since
$a_i(t) <_Q a'_i(t)$, it follows that $\nu(u) < \eta(u)$. As $p$ and $q$
preserve the $P_{\eta}$,
$p(a_i)(t) <_Q q(a'_i)(t)$.
\end{itemize}

To show the generic object is a map defined on all of $A$, it
suffices to show that for any $p\in \Pscr$ and any $a \in A -\dom p$
there is a $q\in \Pscr$ with $p\subseteq q$ and $\dom q = \dom p \union
\{a\}$.
Let $\langle a_1 \ldots a_n \rangle$ enumerate
$\dom p$
in lexicographic order.
Fix $s < n$ with $a_s < a < a_{s+1}$ (the other cases are similar).
Let $m$ be least such that $a_s|(m+1), a|(m+1), a_{s+1}|(m+1)$ are
distinct and let $c$ denote $a|m$.  Suppose
$P_{\rho}(a_s)$,
$P_{\sigma}(a)$, and
$P_{\tau}(a_{s+1})$.
Note that since $a_s < a < a_{s+1}$,
it is impossible for $a_s$ and $a_{s+1}$ to agree on a larger initial
segment
than $a$ and $a_s$ do.  Thus without loss of generality we may assume
that $a_s|m = a|m$.  Two  cases remain.
\begin{itemize}
\item{Case 1:}
$a_s|m = a|m = a_{s+1}|m=c$. Suppose $m$ is odd.
Let $b_s = p(a_s)$
and $b_{s+1} = p(a_{s+1})$.
Then $b_s|m = b_{s+1}|m$ and
$b_s(m) <
b_{s+1}(m)$.
By the definition of amenability for any $r$ with
$b_s(m) < r <
b_{s+1}(m)$, there is an $\eta \in B \inter S_{\sigma}$ with
$\eta|(m+1) = b_s\con r$ so $p\union \{<a,\eta>\} \in \Pscr$
as required.

If $m$ is even choose as the image of $a$, $p(c) \con \sigma(m/2)$.
\item{Case 2:}
$a_s|m = a|m= c$ but $a_{s+1}|m\neq c$.
Again, let
$b_s = p(a_s)$
and $b_{s+1} = p(a_{s+1})$ and denote $b_s|m$ by $b'$.  By amenability
there is an $\eta \in B\inter S_{\sigma}$ with $\eta|m = b'$.  Any such
$\eta$ is less than $b_{s+1}$.  If $m$ is even $\eta > b_s$ is
guaranteed by $\eta(m) > \rho(m)$; if $m$ is odd by
Remark~\ref{recastamen} which
recasts the definition of amenability
we can choose $\eta(m) > b_s(m)$.  In
either case $\eta$ is required image of $a$.
\end{itemize}

We deduce three results from this lemma.  First we note that there
are nonisomorphic but potentially isomorphic suborderings of the reals.
Then we will show in two stages that any countable theory that is not
classifiable has a pair of
models of power $2^{\aleph_0}$
that are not isomorphic but are potentially
isomorphic.

\begin{thm} \label{potiso}
Any two suborderings of $\langle C,<\rangle$ that induce amenable
$L^t$-structures
are potentially isomorphic.
\end{thm}

\Proof.
Since the isomorphism we constructed in proving Lemma \ref{ml}
preserves levels, restricting it to the infinite sequences and reducting
to $<$ yields the required isomorphism.

\jtbdef
\label{treeindis}
Let $M$ be an $L$-structure.
We say that
$\langle \abar_{\eta}\in M:\eta \in Q^{\leq \omega}\rangle$
is a set of {\em $L$-tree indiscernibles}
if for any two sequences $\etabar,
\nubar$ from $Q^{\leq \omega}$:

If $\etabar$ and $\nubar$ realize the
same atomic type in $\langle Q^{\leq\omega};\sqsubset,<,
lev,\meet\rangle$ then
$\langle \abar_{\eta_1}, \ldots \abar_{\eta_n}\rangle$
and
$\langle \abar_{\nu_1}, \ldots \abar_{\nu_n}\rangle$
satisfy the same $L$-type.

\bigskip
Note that the isomorphism given in Theorem~\ref{potiso} preserves
$\wedge$ since $\wedge$ is definable from $\sqsubset$.  We have
introduced $\wedge$ to the language so that atomic types suffice in
Definition~\ref{treeindis}.

\begin{thm}
\label{unss}
Let $T$ be a complete unsuperstable theory in a language
$L$.  Suppose $L \subseteq L_1$ and  $T
\subseteq T_1$ with $|T_1| \leq 2^{\omega}$.  Then there are
$L_1$-structures $M_1, M_2\sat T_1$ such that each
$M_i|L$ is a model of $T$
of cardinality $2^{\omega}$, $M_1$ and $M_2$ are not $L$-isomorphic
but in a ccc-forcing extension of the universe $M_1 \iso_{L_1} M_2$.
\end{thm}
\Proof.  We may assume that $T_1$ is Skolemized.
Note there is no assumption that $T_1$ is stable.
Let $M$ be a
reasonably saturated model of $T_1$. By VII.3.5(2)
of \cite{Shelahbook}
there are $L$-formulas  $\phi_i(\xbar,\ybar)$
for $i \in \omega$ and  a tree of elements
$\langle \abar_{\eta}\in M:\eta \in Q^{\leq \omega}\rangle$ such that
for any $n\in \omega$,
$\eta \in Q^{\omega}$, and
$\nu \in Q^{n+1}$
if $\nu|n = \eta|n$ then
$\phi_{n+1}(\abar_{\eta},\abar_{\nu})
$ if and only if $\nu \sqsubset \eta$.
By
VII.3.6(3)
of \cite{Shelahbook}  (applied in $L_1$!) we may assume that
the index set is a collection of $L_1$-tree indiscernibles.

Let $Y =
\langle \abar_{\nu}\in M:\nu \in Q^{< \omega}\rangle$.  For $\eta \in
Q^{\omega}$, let $p_{\eta}$ be the type over $Y$
containing
$\phi_{n+1}(\xbar;\abar_{\eta|n\con \eta(n)})
\wedge \neg\phi_{n+1}(\xbar;\abar_{\eta|n\con \eta(n)+1})$
for all $n\in \omega$.

Now a direct calculation from the definition of tree indiscernibility
(which was implicit in the proof of Theorem VIII.2.6 of
\cite{Shelahbook})
shows:

{\bf Claim.}
For any $\etabar \in Q^{\omega}$ and any Skolem term $f$, if
$f(\abar_{\eta_1}, \ldots,\abar_{\eta_n})$ realizes $p_\nu$ then some $\eta_i =
\nu$.

Let $M_2$ be the Skolem hull of
$C' = Y \union \{\abar_{\eta}:\eta \in C\}$
where $C$ is chosen as in \ref{maindef}.
Since $Y$ is countable there are at most $2^{\aleph_0}$ embeddings
of $Y$ into $M_2$;
let $f_{\eta}$ for $\eta \in Q^{\omega}$ enumerate them.
For $\eta \in Q^{\omega}$, define $b_{\eta}\in S_{\eta}$
by $b_{\eta}(2n) = \eta(n)$ and $b_{\eta}(2n+1) = 0$ for all $n \in
\omega$.

Let $A = \union_{\eta\in Q^{\omega}}S'_{\eta}$ where
$S'_\eta$ =
$S_{\eta}- \{b_{\eta}\}$ if $M_2$ realizes $f_{\eta}(p_{b_{\eta}})$
and
$S'_\eta$ =
$S_{\eta}$ if $M_2$ omits $f_{\eta}(p_{b_{\eta}})$.

It is easy to check that $A^* \subseteq C^*$ is amenable.  Let $M_1$ be
the Skolem Hull of
$A' = Y \union \{\abar_{\eta}:\eta \in A\}$.

Since $A^*$ and $C^*$ are amenable there is a ccc-forcing notion $\Pscr$
such that $V[G] \sat A^* \iso C^*$.
Since $A'$ and $C'$ are sets of $L_1$-tree indiscernibles, the induced
map is an $L_1$-isomorphism.  Thus, $V[G] \sat M_1
\iso_{L_1} M_2$.  Thus we need only show that $M_1$ and $M_2$ are not
isomorphic in the ground universe.  Suppose $h$ were such an
isomorphism.  Choose $\eta \in Q^{\omega}$ such that $h|Y = f_{\eta}$.
Now if $b_{\eta}\in A$ the construction of $A$ guarantees that $M_2$
omits
$f_\eta(p_{b_{\eta}})
=h(p_{b_{\eta}})$ but $b_{\eta}$ realizes $p_{b_{\eta}} \in M_1$.
On the other hand, if $b_{\eta} \not \in A$ then by the Claim,
$M_1$ omits $p_{b_{\eta}}$ but $M_2$ realizes $f(p_{b_{\eta}})$.
\bigskip

We now want to show the same result for theories with DOP or OTOP.  We
introduce some specialized notation to clarify the functioning of DOP.

\jtbnot
For a structure $M$ elementarily embedded in a
sufficiently saturated structure $M^*$, $\bbar$ from M,
and $\abar$ from $M^*$,
$\dim(\abar,\bbar,M)$ is the minimal cardinality of a
maximal, independent over $\bbar$,
set of realizations of $\stp(\abar/\bbar)$ in $M$.
For  models $M$ of superstable theories, if $\dim(\abar,\bbar,M)$
is infinite
then it is equal to the cardinality of {\it any\/} such maximal
set.
For $p(\xbar,\ybar)
\in S(\emptyset)$ and
$\bbar$ from $M$ let
$d(p(\xbar;\bbar);M)=\sup\{\dim(\abar',\bbar,M):{\rm tp}(\abar'/\bbar)=
p\}$.

\begin{lemm}
\label{dop}
If a complete, first order, superstable theory $T$ of cardinality $\lambda$
has DOP then there is a type $p(\vbar,\ubar,\xbar,\ybar)$ such that
for any cardinal $\kappa$ there is a model $M$ and a sequence
$\{\abar_\alpha:\alpha\in\kappa\}$ from $M$ such that
for all $\alpha,\beta\in\kappa$ and all $\cbar$ from $M$,
$d(p(\vbar;\cbar,\abar_\alpha,\abar_\beta);M)\leq\lambda^+$ and
\begin{equation}
(\exists \ubar \in M)\,
[d(p(\vbar;\ubar,\abar_\alpha,\abar_\beta);M)=
\lambda^+]\  \text{if and only if}\
\alpha<\beta.
\label{onlyeq}
\end{equation}
\end{lemm}

\proof.  This is the content of condition (st 1) on page 517
of \cite{Shelahbook2nd}.
 (As for any infinite indiscernible ${\bf I}$ there is a finite
${\bf J}\subseteq {\bf I}$ such that if
$\dbar\in {\bf I}\setminus {\bf J}$ then
tp$(\dbar, \cup {\bf J}$) is a stationary type
and Av$({\bf I}, \cup {\bf I})$ is a non-forking extension of it).

\begin{prop}
\label{surrogate}
Suppose $|L|=\lambda$ and $T$ is a superstable $L$-theory with either
DOP or OTOP.  There is
an expansion $T_1\supseteq T$,
$|T_1|=\lambda^+$ such that $T_1$ is Skolemized, and
an $L$-type $p$
($p=p(\vbar, \ubar, \xbar, \ybar)$ if $T$ has DOP, $p=p(\vbar,\xbar,\ybar)$
if $T$ has OTOP)
such that $\vbar, \ubar, \xbar, \ybar$ are finite, $\lg \xbar=\lg\ybar$
and
for any order type
$(I,<)$ there is a model $M_I$ of $T_1$ and a sequence $\{\abar_i:i\in I\}$
from $M_I$ of $L_1$-order indiscernibles
such that:

 (a) $M_I$ is the Skolem Hull of $\{\abar_i: i\in I\}$;

(b) If $T$ has DOP then for all $i,j\in I$,
$$(\exists \ubar\in M_I)\,
[d(p(\vbar;\ubar, \abar_i,\abar_j);M_I)\ge
\lambda^+]\  \text{if and only if}\ i<_I j;$$

(c) If $T$ has OTOP then
for all $i, j\in I$,
$M_I\models (\exists \vbar)p(\vbar, \abar_i, \abar_j)$ iff $i<_I j$.
\end{prop}

\proof.
Let $\kappa$ be the Hanf number for omitting types for first-order
languages of cardinality $\lambda^+$.
 If $T$ has OTOP then by its definition (see \cite{Shelahbook2nd} XII \S 4)
there is a model $M$ of
$T$ and  sequence $\{\abar_\alpha:\alpha\in\kappa\}$
of finite tuples
from $M$ and type
$p(\vbar, \xbar, \ybar)$ such that
$M\models (\exists \vbar)
p(\vbar, \abar_\alpha, \abar_\beta)$ iff
$\alpha<\beta$.

  By Lemma~\ref{dop} when $T$ has DOP
 we can find
a model $M$ of $T$, a sequence $\{\abar_\alpha:\alpha\in\kappa\}$ and
a type $p(\vbar,\ubar,\xbar,\ybar)$ so that $(\exists \ubar\in M)\,
[d(p(\vbar;\ubar,\abar_\alpha,\abar_\beta);
M) \ge \lambda^+]$ if and only if $\alpha<\beta$.
We may also
assume that $d(p(\vbar;\cbar,\abar_\alpha,\abar_\beta);M)\leq\lambda^+$
for all $\alpha,\beta\in\kappa$ and $\cbar$.

Let $L_0$ be a minimal Skolem expansion of $L$.  That is,
$L_0$ is a minimal expansion of $L$ such that
there is a function
symbol $F_\phi(\ybar)\in L_0$ for each formula $\phi(x,\ybar)\in L_0$.
Let $M_0$ be any expansion of $M$ satisfying the Skolem axioms
$\forall\ybar [(\exists x)\phi(x,\ybar)\rightarrow\phi(F_\phi(\ybar),\ybar)]$
and let $T_0={\rm Th}(M_0)$. Without loss
of generality $\lambda^++1\subseteq M_0$.

>From now on, assume we are in the DOP case as the OTOP case is
similar and does not require a further expansion of the language
(i.e., take $L_1=L_0$ and $T_1=T_0$.)
Expand $L_0$ to $L_0'$ by adding relation symbols $<, \in, P,$
constants for all ordinals less than or equal to $\lambda^+$
and
a new function symbol
$f(w,\ubar,\xbar,\ybar)$.
Let $M_0'$ be an expansion of
$M_0$ so that $<$ linearly orders the $\abar_{\alpha}$ and
the set of $\abar_{\alpha}$ is the denotation of $P$.
Interpret the constants and $\in$ in the natural way.
For all $\alpha,\beta\in\kappa$ and all realizations $\dbar\conc\cbar$ of
$p(\vbar,\ubar, \abar_\alpha,\abar_\beta)$ in $M_0$,
let $(\lambda w)f(w,\cbar,\abar_{\alpha},\abar_{\beta})$
be a 1--1 map from an initial segment of $\lambda^+$
to a maximal,
independent over $\cbar\cup\abar_\alpha\cup\abar_\beta$, set of realizations of
$\stp(\dbar/\cbar\,\abar_{\alpha}\abar_{\beta})$.

Let $L_1$ be a minimal Skolem
expansion of $L_0'$, let
$M_1$ be a Skolem expansion of
$M_0'$ to an $L_1$-structure and let
$T_1$ denote the theory of $M_1$.  So $|T_1|=\lambda^+$.

Note that if, for some $\cbar$, the domain of $(\lambda w)
f(\cbar,\abar_\alpha,\abar_\beta)$ is $\lambda^+$ then $\alpha <\beta$.
Also, for all $\alpha,\beta\in\kappa$ and $\cbar$ from $M_1$
the independence of the range of
$(\lambda w)
f(w,\cbar,\abar_{\alpha},\abar_{\beta})$ is expressed by
an $L_1$-type.
Thus $M_1$ omits the types
$$q(\vbar; \ubar,\xbar,\ybar) = p(\vbar,\ubar, \xbar,\ybar)
\union \{\fg{\vbar}{\ubar\xbar\ybar}{\{
f(\gamma,\ubar,\xbar,\ybar)
:\gamma < \lambda^+ \}}
\} \union \{P(\xbar)\} \union \{P(\ybar)\}\cup\{\xbar\not <\ybar\}
$$ and $r(v) = \{v \in \lambda^+\} \union \{v \neq \gamma:\gamma
\in \lambda^+\}$.

To complete the proof of the proposition construct
an Ehrenfeucht-Mostowski model $M_I$ of $T_1$
built from a set of $L_1$-order indiscernibles $\{\abar_i:i\in I\}$
omitting both  $q(\vbar,\ubar,\xbar,\ybar)$ and $r(\vbar)$.
The existence of such a model follows as in
the proof of Morley's
omitting types theorem (see e.g.,
\cite[VII.5.4]{Shelahbook}).

\bigskip
Note that in the DOP case of the proposition above
the argument shows
$d(p(\vbar;\cbar,\abar_i,\abar_j);M_I)\le\lambda^+$ for all $i,j\in I$ and
$\cbar$.

We have included a sketch of the proof of Lemma~\ref{surrogate} which is
essentially Fact X.2.5B+$620_9$ of \cite{Shelahbook2nd} and
Theorem~0.2 of \cite{Shelah220}  to clarify two points.
We would not include
this had not experience showed that some readers miss these points.
Note that the parameter $\cbar$ is needed in the DOP case
not only to fix the strong type, but because in general we cannot ensure the
existence of a large, independent set of realizations over $\abar_\alpha
\cup\abar_\beta$.
  Also, it is essential
that we pass to a Skolemized
expansion to carry out the omitting types argument
and that the final
set of indiscernibles are indiscernible in the Skolem language.
We can then reduct
to $L$ for the many models argument
(if we use
III 3.10 of \cite{Shelahuniversal}) not just \cite{Shelahbook2nd} VIII,\S3)
 but for the purposes of
this paper we cannot afford to take reducts as the proof of Theorem~1.14
requires that an isomorphism between linear orders $I_1, I_2$
induces an isomorphism of the corresponding models.

Let us expand on why we quote \cite{Shelahuniversal} above.
In \cite{Shelahuniversal},~Theorem~III~3.10 it is proved that for
all uncountable cardinals $\lambda$ and all vocabularies $\tau$,
if there is a formula $\Phi(\xbar,\ybar)$ such that for every
linear order $(J,<)$ of cardinality $\lambda$ there is a $\tau$-structure
$M_J$ of cardinality $\lambda$
and a subset of elements $\{\abar_s:s\in J\}$
satisfying
\begin{enumerate}
\item  $M_J\models \Phi(\abar_s,\abar_t)$ if and only if $s<_J t$ and
\item  The sequence $\langle\abar_s:s\in J\rangle$ is skeleton-like in
$M_J$ (i.e., any formula of the form $\Phi(\xbar,\bbar)$ or $\Phi(\bbar,\xbar)$
divides
$\langle\abar_s:s\in J\rangle$ into finitely many intervals)
\end{enumerate}
then there are $2^\lambda$ non-isomorphic $M_J$'s.

The point, compared with earlier many-models
proofs, is that we do not demand that
the $M_J$'s be constructed from $J$ in any specified way.  It is true
that the natural example
satisfying these conditions is an Ehrenfeucht-Mostowski model
built from $\langle\abar_s:s\in J\rangle$
in some expanded language, but this is not required.  In particular,
our generality allows taking reducts, so long as the formula $\Phi$
remains in the vocabulary.  Further, there is no requirement that
$\Phi$ be first-order.

However, in our context we want to introduce an isomorphism between
two previously non-isomorphic models.  The natural way of doing this
is to produce two non-isomorphic but potentially isomorphic orderings
$J_1$ and $J_2$ and then conclude that $M_{J_1}$ and $M_{J_2}$ become
isomorphic.  Consequently, it is important for us to know that the models
are E.M.\ models.

\bigskip

We can simplify the statement of the conclusion of
Lemma~\ref{surrogate} if we define the logic
with `dimension quantifiers'.
In this logic we demand that in addition to the requirement that
`equality' is a special predicate to be interpreted as identity that
another family of predicates also be given a canonical interpretation.

\jtbnot
Expand the vocabulary $L$ to $\hat L$ by
adding
new predicate symbols $Q_\mu(\xbar,\ybar)$ of each finite arity
for all cardinals $\mu\le\lambda^+$.
Now define the logic
$\hat{L}_{\lambda^+,\omega}$ by first demanding that each predicate
$Q_\mu$ is interpreted
in an $L$-structure
$M$ by $$M\models Q_\mu(\abar,\bbar)\ \text{if and only if}\
\dim(\abar,\bbar,M)=\mu.$$
Then define the quantifiers and connectives as usual.
We will only be concerned with the satisfaction of sentences of this
logic for models of superstable theories.

\jtbnumpar{Remarks}
\label{real}
\begin{enumerate}
\item \label{absolute}
The property coded in
Condition~(\ref{onlyeq}) of Lemma~\ref{surrogate} is expressible
by a formula $\Phi(\xbar,\ybar)$ in the logic
$\hat{L}_{\lambda^+,\omega}$.
Each formula in $\hat L_{\lambda^+,\omega}$ and in particular
this formula $\Phi$ is
absolute relative to
any extension of the universe that preserves cardinals.
More precisely
$\Phi$ is
absolute relative to
any extension of the universe that preserves
$\lambda^+$.
\item
  If $T$ has OTOP the formula $\Phi$ can be taken in the logic
$L_{\lambda^+,\omega}$.
So in this case $\Phi$ is preserved in any forcing extension.
\item
Alternatively, the property coded in Condition (1) of
Lemma~\ref{surrogate} is also expressible in $L_{\lambda^+,\lambda^+}$.
That is, there is a formula $\Psi(\xbar,\ybar)\in L_{\lambda^+,\lambda^+}$
(in the original vocabulary $L$) so that
$$M_I\models \Psi(\abar_i,\abar_j)\ \text{if and only if}\ i<_I j.$$
The reader should note that satisfaction of arbitrary sentences of
$L_{\lambda^+,\lambda^+}$ is, in general, not absolute for
cardinal-preserving forcings.
However, the particular statements $M_I\models \Psi(\abar_i,
\abar_j)$ and $M_I\models \neg\Psi(\abar_i,\abar_j)$ will be preserved
under any cardinal-preserving forcing by the first remark.
\item
Note that we could have chosen the type $p$
(in the DOP case) such that
$p(\vbar; \cbar, \abar_\alpha, \abar_\beta)$ is a stationary regular type.
Note also that had we followed \cite{Shelahbook2nd},X2.5B more closely,
we could have insisted that $\vert T_1\vert=\lambda$.
In fact we could have arranged that in $M_I$,
every dimension would be $\le \aleph_0$
or
$\Vert M_I\Vert$
(over a countable set).
However, neither of these observations improve the statement of \ref{nonclass}.

\end{enumerate}

\begin{thm}
\label{nonclass}
If $T$ is a complete theory in a vocabulary $L$
with $|L|\leq 2^{\omega}$ and
$T$ has either OTOP or DOP then there are models $M_1$ and $M_2$ of $T$
with
cardinality the continuum that are not isomorphic but are potentially
isomorphic.

\end{thm}

\Proof.  By Theorem~\ref{unss} we may assume that $T$ is superstable.
By Proposition~\ref{surrogate} and Remark~\ref{real}(\ref{absolute})
there is a model
$M$ of a theory $T_1\supseteq T$ in a
Skolemized language $L_1\supseteq L$
containing a set of $L_1$-order indiscernibles
$\{\abar_{\eta}:\eta \in Q^{\leq\omega} \}$
and an $\hat{L}_{\lambda^+,\omega}$-formula $\Phi(\xbar,\ybar)$
so that $\Phi(\abar_{\eta},\abar_{\nu})$ holds in $M$
if and only if  $\eta$ is
lexicographically less than $\nu$.
Further, the statements ``$M\models \Phi(\abar_\eta,\abar_\nu)$''
and ``$M\models \neg\Phi(\abar_\eta,\abar_\nu)$'' are preserved under
any ccc forcing.
Note that this $L_1$-order indiscernibility certainly
implies $L_1$-tree-indiscernibility in the sense of
Definition~\ref{treeindis}.

Thus, the construction of potentially isomorphic but not isomorphic models
proceeds as in the last few paragraphs of
the proof of Theorem~\ref{unss} once we
establish the following claim.

\smallskip
{\bf Claim.}  For any $\nu \in Q^{\omega}$ there is a collection
$p_{\nu}(\xbar)$ of boolean combinations of $\Phi(\xbar,\abar)$ as
$\abar$ ranges over $Y$ such that
for any $\etabar \in Q^{\omega}$ and any $L_1$-term $f$, if
$f(\abar_{\eta_1}, \ldots \abar_{\eta_n})$ realizes $p_\nu$ in $M$
then some $\eta_i =
\nu$.

\smallskip
\Proof.
The conjunction of the
$\Phi(\xbar;\abar_{\nu|n \con
<\nu(n) +1>})$    and
$\neg\Phi(\xbar;\abar_{\nu|n \con
<\nu(n)+1>})$ that define the `cut' of $\abar_{\nu}$ will constitute
$p_{\nu}$. Now if $\nu$ is not among the $\eta_i$ choose any $n$ such
that $\eta_1|n, \eta_2|n, \ldots \eta_k|n, \nu|n$ are distinct.
Then the sequences
$\langle \eta_1, \ldots \eta_k,
\nu|n\con\langle\nu(n) +1\rangle\rangle$
and $\langle \eta_1, \ldots \eta_k, \nu|n\con\langle\nu(n) -1\rangle\rangle$
have the same type in the lexicographic order so
$$M\sat \Phi(f(\abar_{\eta_1}, \ldots \abar_{\eta_k});\abar_
{\nu|n\con<\nu(n) +1>})
\iff \Phi(
f(\abar_{\eta_1}, \ldots \abar_{\eta_k});\abar_
{\nu|n\con<\nu(n) -1>}).$$
Thus, $f(\abar_{\eta_1}, \ldots \abar_{\eta_k})$ cannot realize
$p_{\nu}$.

\jtbnumpar{Remarks}
\begin{enumerate}
\item
Note that in Theorem~\ref{unss} we were able to use
any expansion of $T$ as $T_1$ so the result is actually for
$PC_{\Delta}$-classes.  In Theorem~\ref{nonclass} our choice
of $T_1$ was constrained, so
the result is true for only elementary as opposed to
pseudoelementary classes.  The case of unstable elementary classes could
be handled by the second method thus simplifying the combinatorics at
the cost of weakening the result.
\item While we have dealt only with models and theories
of cardinality $2^{\omega}$,
the result extends immediately to models of any larger cardinality
and
straightforwardly to theories of cardinality
$\kappa$ with
$\kappa^{\omega} = \kappa$.
\end{enumerate}

 \section{Classifiable examples}

We begin by giving an example of a classifiable theory
having a pair of non-isomorphic, potentially isomorphic models.
We then extend this result to a class of
weakly minimal theories.

Let the language $L_0$ consist of a countable family $E_i$ of binary
relation symbols and let the language $L_1$ contain an additional
uncountable set of unary predicates
$P_{\eta}$.
We first construct an $L_0$-structure that is rigid but can be
forced by a ccc-forcing to be nonrigid.
Our example will be in the language $L_0$ but we will
use expansions of the $L_0$-structures to $L_1$-structures in the
argument.

We now revise the definitions leading up to the notion of an amenable
structure in Section~\ref{unclass} by replacing the underlying structure
on $Q^{\leq\omega}$ by one
with universe $2^{\omega}$.   In particular, $D_{\eta}$, $S_{\eta}$, and
$C$
are now being redefined.

\jtbnot \mbox{}
\label{maindef2}
\begin{enumerate}
\item
For $\eta \in 2^{\omega}$, let
$D_{\eta} = \{\sigma \in 2^{\omega}:\sigma(2n) =
\eta(n)\}$ and
$S_{\eta} = \{\sigma \in D_{\eta}:\sigma(2n+1)
\text{ is 0 for all but finitely many $n$}\}\cup\{b_\eta\}$,
where $b_\eta$ is any
element of $D_\eta$ satisfying $b_\eta(2n+1)=1$ for infinitely many $n$.
Let $C = \bigunion_{\eta\in 2^{\omega}}S_{\eta}$.
\item
Let $M^*$ be the
$L_1$-structure with universe $2^{\omega}$ where
$E_i(\sigma,\tau)$ holds if $\sigma|i = \tau|i$, and the unary relation
symbol $P_{\eta}$ holds of the set $S_{\eta}$.
Let $M_1$ be the $L_1$-substructure of $M^*$ with universe $C$.
\item  Any subset $A$ of $C$ inherits a natural $L_1$ structure from
$M_1$
with $P_{\eta}$ interpreted as $S_{\eta}\inter A$.
\end{enumerate}

\jtbdef  An $L_1$-substructure $M_0$ of $M_1$
is {\em amenable} if for all $\eta \in 2^{\omega}$, all $n\in \omega$
and all $s \in 2^n$, if there is a $\nu \in P_{\eta}(M_1)$ with $\nu|n =
s$ then there is a $\nu'\in P_{\eta}(M_0)$ with $\nu'|n = s$.

Note that any $L_1$-elementary substructure of $M_1$
is amenable.
Moreover, it easy to see that i)
each $D_{\eta}$ is a perfect tree,
ii) $2^{\omega}$
is a disjoint union of the $D_{\eta}$ and iii)
for each $s\in 2^{<\omega}$ there are
$2^{\omega}$
sequences $\eta$ such that $s$ has an extension $b \in D_{\eta}$.

\begin{thm}
\label{examp}
The theory $FER_{\omega}$ of countably many refining
equivalence relations with binary splitting has a pair of models of size
the continuum which are not isomorphic but are potentially isomorphic.
\end{thm}

This result follows from the next two propositions and the fact that
$M_1$ is not
rigid.
\begin{prop}  There is an $L_1$-elementary
substructure $M_0$ of $M_1$ such that
\begin{enumerate}
\item $|P_{\eta}(M_1) - P_{\eta}(M_0)| \leq 1$.
\item
$M_0|L_0$ is rigid.
\end{enumerate}
\end{prop}

\Proof.
Note that each automorphism of $M_1|L$ is determined by its restriction
to the eventually constant sequences so there are only $2^{\omega}$
such. Thus we may
let $\langle
f_i:i<2^{\omega}\rangle$ enumerate the nontrivial
automorphisms of
$M_1|L$.
We define by induction disjoint subsets $A_i, B_i$ of $M_1$ each with
cardinality
less than the continuum.  We denote $\union_{i<j} A_i$ by $\Aun_j$.
At stage $i$, choose $\alpha \in M_1$ such that $f_i$ moves $\alpha$.
Then, by continuity,
there is a
finite sequence $s$ such that every element of $W_s
=\{\tau:s\subseteq \tau\}$ is moved by $f_i$.
Since $|A_i|, |B_i| < 2^{\omega}$ and by condition iii) of
the remark after the definition of amenable there
are an $\eta
\in 2^{\omega}$ and a $\beta \in P_{\eta} \inter (f_i(W_s) - \Aun_i)$.
Then let
$B_i = \{\beta\}$ and $A_i = S_{\eta}-\{\beta\}$.
Finally,
let $M_0 = M_1 - \Bun_{2^{\omega}}$.

Since no element is ever removed from an $A_i$, condition i)
is satisfied.
It is easy to see that $M_0$ is rigid, as any nontrivial automorphism
$h$ of $M_0$
would extend in a unique way to an automorphism $f_i$
of $M_1$ but at step $i$
we ensured that the restriction of $f_i$ to $M_0$ is not an
automorphism.

\begin{prop}  If $M_0$ is an amenable substructure of $M_1$, $M_0$
and $M_1$ are potentially isomorphic.
\end{prop}

\Proof.  Let $\Pscr$ be the collection of all finite partial
$L_1$-isomorphisms between $M_0$ and $M_1$.

We first claim that $\Pscr$ is a ccc set of forcing conditions.  In
fact, $\Pscr = \cup_{n\in \omega} F_n$ where if $p,q \in F_n$ then $p
\union q \in \Pscr$.  Given $p \in \Pscr$, fix an (arbitrary)
enumeration
$\langle a_1 \ldots a_n \rangle$ of $\dom p$.  Let $n(p)$ be
the cardinality of $\dom p$ and let $k(p)$ be the least integer
satisfying the following conditions.

\begin{enumerate}
\item If $i \neq j$, $M_0 \sat \neg E_k(a_i,a_j)$.
\item If $i \neq j$, $M_1 \sat \neg E_k(p(a_i),p(a_j))$.
\item For all $n$ with $2n +1 \geq k$, $a_i(2n+1) = 0$.
\item For all $n$ with $2n +1 \geq k$, $p(a_i)(2n+1) = 0$.
\end{enumerate}

Now define an equivalence relation on $\Pscr$ by $p \sim q$ if $n(p) =
n(q)$, $k(p) = k(q)$ and
letting
$\langle a_1 \ldots a_n \rangle$ enumerate
$\dom p$,
$\langle a'_1 \ldots a'_n \rangle$ enumerate
$\dom q$, for each $i$,
$a_i|k(p) = a'_i|k(p)$
and
$p(a_i)|k(p) = q(a'_i)|k(p)$.  Then $\sim$ has only countably many
equivalence classes; these classes are the $F_n$.  We must show
that if $p \sim q$ then $p \union q \in \Pscr$.  It suffices to
show that for all $i,j$ if
$M_0 \sat E_n(a_i,a'_j)$ then
$M_1 \sat E_n(p(a_i),q(a'_j))$.
\begin{itemize}
\item{Case 1: $i\neq j$.}  By the definition of $k = k(p)$,
$M_0 \sat \neg E_k(a_i,a'_j)$.
Let $\ell$ be maximal so that
$M_0 \sat E_{\ell}(a_i,a'_j)$.
Then $\ell \geq n$ since $M_0 \sat E_n(a_i,a'_j)$.
Since $a_j|k = a'_j|k$, $\ell$ is also maximal so that $M_0
\sat E_{\ell}(a_j,a_i)$.
As $p$ is an isomorphism,
$M_1 \sat E_{\ell}(p(a_i),p(a'_j))$.
But $p(a_j)|k = q(a'_j)|k$, so $\ell$ is also maximal with
$M_1 \sat E_{\ell}(q(a'_j),p(a_i))$.  Since $n \leq \ell$
we conclude $M_1 \sat E_n(q(a'_j),p(a_i))$ as required.
\item{Case 2: $i = j$.}  Suppose $a_i \in S_{\nu}$ and $a'_j \in
S_{\eta}$.  We have $M_0 \sat E_k(a_i,a'_j)$ by the definition of $k$
and similarly, $M_1 \sat E_k(p(a_i),q(a'_i))$.  Now we show by induction
that for
each $m \geq k$, $M_0 \sat E_m(a_i,a'_i)$ if and only if
$M_1 \sat E_m(p(a_i),q(a'_i))$.
Assuming this condition for $m$ we show it for $m+1$.  If $m+1$ is odd, the
result is immediate by parts iii) and iv) of the conditions defining
$\sim$.  If $m = 2u$ then
$a_i(m) =\nu(u)$ and $a'_i(m) = \eta(u)$.
Since $p$ and $q$ are $L_1$-isomorphisms
$p(a_i)(m) =\nu(u)$ and $q(a'_i)(m) = \eta(u)$.  But
$M_0 \sat
E_m(a_i, a'_i)$
implies $\nu(u) = \eta(u)$ so we have
$M_1 \sat E_m(p(a_i), q(a'_i))$.
\end{itemize}

To show that the generic object is a map with domain $M_0$, it
suffices to show that for any $p\in \Pscr$ and any $a \in M_0 -\dom p$
there is a $q\in \Pscr$ with $q\leq p$ and $\dom q = \dom p \union
\{a\}$.  Choose $n$ so that the members of $\dom p \union \{a\}$ are
pairwise $E_n$-inequivalent.  Fix any $L_1$-automorphism $g$ of
$M^*$ that extends $p$.
Let $s = g(a)|n$ and $W_s =\{\gamma \in 2^{\omega}:s \subseteq
\gamma\}$.
Since there is a $\nu \in C \inter W_s$ and $B$ is amenable, there is a
$\nu' \in B \inter W_s$.
Choosing $\nu'$ for $b$, $p\union \langle
a,b\rangle$ is the required extension of $p$.

Since $M_0$
and $M_1$ are isomorphic in a generic extension for this forcing, we
complete the proof.
\bigskip

\jtbnumpar{Remark}  The notion of a classifiable theory having two
non-isomorphic, potentially isomorphic models is not very robust,
and in particular can be lost by adding constants.
As an example, let $FER^*_\omega$ be an expansion of $FER_\omega$
formed by adding constants for the elements of a given countable model of
$FER_\omega$.  Then every type in this expanded language is stationary
and the isomorphism type of any model of $FER_\omega^*$ is determined
by the number of realizations of each of the $2^\omega$ non-algebraic
1-types.  Thus, if two models of $FER^*_\omega$ are non-isomorphic then
they remain non-isomorphic under any cardinal-preserving forcing.

Similarly, non-isomorphism of models of the theory $CEF_\omega$
of countably many crosscutting equivalence relations (i.e.,
Th$(2^\omega,E_i)_{i\in\omega}$, where $E_i(\sigma,\tau)$ iff
$\sigma(i)=\tau(i)$) is preserved under ccc forcings.

\bigskip
We next want
to extend the result from Theorem~\ref{examp}  to a larger class of
theories.  Suppose $T$ is superstable and there is a type $q$, possibly
over a finite set $\ebar$
of parameters, and an $\ebar$-definable family $\{E_n:n\in\omega\}$
of properly refining equivalence relations, each with finitely
many classes that determine the strong types extending $q$.
Let $T$ be such a theory in a language $L$ and let $M$ be a model of $T$.
Let $L_0$ be a reduct of $L$ containing the $E_n$'s.

We say
$\langle a_{\eta}\in M:\eta \in  X\subseteq 2^{\omega}\rangle$
is a set of {\em unordered
tree $L$-indiscernibles}
if the following holds
for any two sequences $\etabar,
\nubar$ from $X$:

If $\etabar$ and $\nubar$ realize the
same $L_0$-type
then
$\langle a_{\eta_1}, \ldots a_{\eta_n}\rangle$
and
$\langle a_{\nu_1}, \ldots a_{\nu_n}\rangle$
satisfy the same $L$-type.

  We say that a superstable theory $T$ with a type of
infinite multiplicity as above {\em embeds an unordered tree}
if there is a model $M$ of $T$ containing a set of unordered
tree $L$-indiscernibles indexed by $2^\omega$.
We deduce below the existence of potentially isomorphic nonisomorphic
models of weakly minimal theories which embed an unordered tree.
Every small superstable, non-$\omega$-stable theory has a type of
infinite multiplicity with an associated family of $\{E_n:n<\omega\}$ of
refining equivalence relations and a set of tree indiscernibles in the
sense of \cite{Baldwindivers}.  The existence of such a tree of
indiscernibles suffices for the many model arguments but does not in
itself suffice for this result.  Marker has constructed an
example of such a theory which does not embed an unordered tree.
However, an apparently ad hoc argument shows this example does have
potentially isomorphic but not isomorphic models.

\jtbnot
\label{treenot}
Given $A=\{a_\eta:\eta\in 2^\omega\}$ a set of unordered
tree indiscernibles
let $D=\{a_\eta\in A:\eta(n)=0$ for
all but finitely many $n\}$.  For $\eta\in 2^\omega$ let
$p_\eta(x)\in S^1(D)$ be $q(x)\cup \{E_n(x,a_\nu):a_\nu\in D$ and $\nu
|n=\eta|n\}$.  Note that $D$ is a dense
subset of $A$, each $a_\eta$ realizes $p_\eta$ and
each $p_\eta$ is stationary.

\begin{lemm}
\label{setstage}
Let $T$ be a weakly minimal theory that embeds an unordered tree.
Fix $A$ and $D$ as described in Notation~\ref{treenot}.
There is a set $X$ satisfying the following conditions:
\begin{enumerate}
\item $X \union A$ is independent over the empty set;
\item
for any $Y$ with $D \subseteq Y \subseteq A$, and
any $\eta \in
2^{\omega}$,
$p_{\eta}$ is realized in $\acl(XY)$ if and only if $p_{\eta}$ is
realized in $Y$;
\item
for any $Y$ with $D \subseteq Y \subseteq A$,
$\acl(XY)$ is a model
of $T$.
\end{enumerate}
\end{lemm}

\proof.  It is easy to see from the definition of unordered tree
indiscernibility that if $X = \emptyset$, then conditions i) and ii)
of the
Lemma are satisfied for any $Y \subseteq A$.
We
will show that for any $X$ and $Y$ with $D\subset Y \subseteq A$ with
$XY$ satisfying conditions i) and ii) and any consistent formula
$\phi(v)$ over $\acl(XY)$
that is not satisfied in $\acl(XY)$
it is possible to adjoin a solution of $\phi$
to $X$ while preserving the conditions. By iterating this procedure we
obtain a model of $T$.

Now suppose
there is a $Y$ with $D \subseteq Y \subseteq A$,
such that
$\acl(XY)$ is not an elementary submodel of the monster.
Choose a formula $\phi(x,\cbar,\abar)$
with $\cbar \in X$ and
$\abar \in Y$ such that $\phi(x,\cbar,\abar)$ has a solution $d$ in
$\Mscr$ but not in $\acl(XY)$.  If we adjoin $d$ to $X$ we must check
that conditions i) and ii) are not violated.
Since $T$ is weakly minimal and $d\not\in
\acl(XY)$, $XAd$ is independent.
Suppose for
contradiction
that for some $\abar' \in Y$,
$p_{\nu}$ is not realized in  $\abar'$ but $p_{\nu}$ is
realized in $\acl (Xd\abar')$ by say $e$.
Since condition ii) holds for $XY$, $e\not\in \acl(X\abar')$.
Therefore by the exchange lemma
$d \in \acl(Xe\abar')$.
Let $\theta(v,\cbar',\abar,'e)$ with $\cbar' \in X$ and $\abar'\in Y$
witness this algebraicity.  Then
$$\chi(\cbar,\cbar',\abar,\abar',z) = (\exists x) [\phi(x,\cbar,\abar)
\wedge \theta(x,\cbar',\abar',z)]
\wedge (\exists ^{=m}x)\theta(x,\cbar',\abar',z)$$
is a formula over
$X\abar\abar'$ satisfied by $e$.
Moreover, $e\not \in \acl(X\abar\abar')$. For, if so, transitivity
would give $d\in \acl(X\abar\abar') \subseteq \acl(XY)$.
Now $\tp(e/X\abar\abar')$ and in particular
$\chi(\cbar,\cbar',\abar,\abar',z)$
is implied by $p_{\nu}$ and the assertion that
$z\not \in \acl(X\abar\abar')$.  Since $XA$ is independent,
it follows by compactness that there is $b \in D$ such that
$\chi(\cbar,\cbar',\abar,\abar',b)$ holds.
So there is a solution
of $\phi(x,\cbar,\abar)$ in the algebraic closure of $XY$.
This contradicts the original choice of $\phi$ so we conclude that
condition ii) cannot be violated.

\begin{thm}  If $T$ is a  weakly minimal theory in a language
of cardinality at most $2^{\aleph_0}$
that embeds an
unordered tree
then $T$ has two models
that are not isomorphic but are potentially isomorphic (by a
ccc-forcing).
\end{thm}

\Proof.  Let $L$ be the language of $T$.
Assume that the type $q$ is based on a finite set $\ebar$.  Let
$T'$ be the expansion of $T$ formed by adding constants for $\ebar$.
Let $\Mscr$ be a large saturated model of the theory
$T'$ and let the
sets $A$ and $D$  be chosen as in
Lemma~\ref{setstage} applied in $L'$ to $T'$.

Recall the definition of $C$ from Notation~\ref{maindef2}. For
any $W\subseteq C$, let $M'_W$ be the $L'$-structure with universe
$\acl(X \union \{a_{\eta}:\eta \in
W\})$
and denote $M'_W|L$ by $M_W$.
We will construct an amenable set
$W$
such
that $M_W \not \iso M_C$.
Since both are amenable, there is a forcing
extension where $W \iso C$ as $L_1$-structures.
Since $\{a_{\eta}:\eta \in C\}$ is a set of unordered tree
$L'$-indiscernibles, the induced mapping of
$\{a_{\eta}:\eta \in W\}$ into
$\{a_{\eta}:\eta \in C\}$ is $L'$-elementary.  Thus,
$M'_W \iso_{L'} M'_C$
and a fortiori
$M_W \iso_{L} M_C$.

To construct $W$, let $\{f_{\eta}:\eta \in 2^{\omega}\}$ enumerate all
$L$-embeddings of $D\ebar$
into $M_C$.
Note that each $p_{\eta}$ can be considered as a complete $L$-type over
$D\ebar$.

Let $W = \union_{\eta\in 2^{\omega}}S'_{\eta}$ where
$S'_\eta$ =
$S_{\eta}- \{b_{\eta}\}$ if $M_C$ realizes $f_{\eta}(p_{b_{\eta}})$
and
$S'_\eta$ =
$S_{\eta}$ if $M_C$ omits $f_{\eta}(p_{b_{\eta}})$.\ (see Notation~2.1.)

Suppose for contradiction that $g$ is an $L$-isomorphism
between $M_W$ and $M_C$.  Then for some $\eta$, $g|D =
f_{\eta}$.  Now if $c_\eta \in W$, the definition of $W$ yields
$f_{\eta}(p_{\eta})$ is not realized in $M_C$.  This contradicts the
choice of $g$ as an isomorphism.  But if $c_{\eta}$ is not in $W$ then
by the construction of $W$, $f_{\eta}(c_\eta) = g(c_{\eta})$ does not
realize $g(p_{\eta})$.  But this is impossible since $g$ is a
homomorphism.

\bibliography{ssgroups}

\begin{thebibliography}{1}

\bibitem{Baldwindivers}
J.T. Baldwin.
\newblock Diverse classes.
\newblock {\em Journal of Symbolic Logic}, 54:875--893, 1989.

\bibitem{BuechlerShelah}
Steve Buechler and Saharon Shelah.
\newblock On the existence of regular types.
\newblock {\em Annals of Pure and Applied Logic}, 45:207--308, 1989.

\bibitem{Hartthesis}
Bradd Hart.
\newblock {\em Some results in classification theory}.
\newblock PhD thesis, McGill University, 1986.

\bibitem{Shelahbook}
S.~Shelah.
\newblock {\em Classification {T}heory and the {N}umber of {N}onisomorphic
  {M}odels}.
\newblock North-Holland, 1978.

\bibitem{Shelah220}
S.~Shelah.
\newblock Existence of many ${L}_{\infty,\lambda}$-equivalent non-isomorphic
  models of ${T}$ of power $\lambda$.
\newblock {\em Annals of Pure and Applied Logic}, 34, 1987.

\bibitem{Shelahuniversal}
S.~Shelah.
\newblock Universal classes: Part 1.
\newblock In J.~Baldwin, editor, {\em Classification Theory, Chicago 1985},
  pages 264--419. Springer-Verlag, 1987.
\newblock Springer Lecture Notes 1292.

\bibitem{Shelahbook2nd}
S.~Shelah.
\newblock {\em Classification {T}heory and the {N}umber of {N}onisomorphic
  {M}odels}.
\newblock North-Holland, 1991.
\newblock second edition.

\end{thebibliography}
\bibliographystyle{plain}
 \end{document}